\documentclass[12pt,a4paper,draft]{article}
\usepackage[english]{babel}
\usepackage{amssymb, amsmath, latexsym, amsfonts, amsthm}
\begin{document}
\newtheorem{theorem}{Theorem}
\newtheorem{lemma}{Lemma}
\begin{center}
ON A SPACE OF RAPIDLY DECREASING 

INFINITELY DIFFERENTIABLE FUNCTIONS 

ON AN UNBOUNDED CONVEX SET IN ${\mathbb R}^n$ AND ITS DUAL 
\footnote
{This work was supported by the grants RFBR 08-01-00779, 08-01-97023 and LSS-3081.2008.1.}
\end{center}
\begin{center}
I.Kh. MUSIN, P.V. YAKOVLEVA
\end{center}
\renewcommand{\abstractname}{}
\begin{abstract}
{\sc Abstract}.  
Description of linear continuous  functionals on a space of rapidly decreasing infinitely differentiable functions on an unbounded closed convex set in $\mathbb R^n$ in terms of their Fourier-Laplace transform is obtained.  
\end{abstract}

\vspace {0,3cm}

MSC-class: 46F05, 32A10. 
\vspace {0,5cm}

\begin{center}
{\bf 1. Introduction} 
\end{center}

{\bf 1.1. On a problem}. 
Let $C$ be an open convex acute cone in ${\mathbb R}^n$ with the apex at the origin [1, p. 73] and  
$b$ be a convex continuous positively homogeneous function of degree 1 on the closure ${\overline C}$ of $C$ in ${\mathbb R}^n$. 
The pair $(b, C)$ determines the closed convex unbounded set
$$
U(b, C) = \{\xi \in {\mathbb R}^n: \ 
 -<\xi, y> \ \le b(y),  
\forall y \in C \},
$$  
not containing a whole line. 
Note that the interior of the set $U(b, C)$ is not empty and coincides with the set 
$$
V(b, C)= \{\xi \in {\mathbb R}^n: \ 
-<\xi, y> \ 
 < b(y),  
\forall y \in {\overline C}\}, 
$$
and closure of $V(b, C)$ in ${\mathbb R}^n$ is $U(b, C)$. 
For brevity denote $U(b, C)$ by $U$ and $V(b, C)$ by $V$. 

Let $M=(M_k)_{k=0}^{\infty} $ be a non-decreasing sequence of numbers $M_k$ such that $M_0=1$ and 
$$
\displaystyle\lim_{k\rightarrow\infty}\frac {\ln M_k} {k} = +\infty.
$$

For $m\in\mathbb N$ and $\varepsilon > 0$ let $G_{m,\varepsilon}(U)$ be the space  of $C^{\infty}(U)$-functions $f$ with a finite norm
$$ 
p_{m, \varepsilon}(f)= \sup_{x \in V, \alpha \in {\mathbb Z_+^n}} \frac {\vert (D^{\alpha}f)(x)\vert (1+\Vert x \Vert)^m}{\varepsilon^{\vert \alpha \vert} M_{\vert \alpha \vert}} \ .
$$ 
Put $G_M(U)=\displaystyle \bigcap_{m=1}^{\infty}\bigcap_{\varepsilon > 0} G_{m, \varepsilon}(U)$. Thus, $G_M(U)$ is a subclass of the Schwartz class $S(U)$ of rapidly decreasing functions on $U$. 

With usual operations of addition and multiplication by complex numbers $G_M(U)$ becomes a linear space. 
The family of norms $p_{m, \varepsilon}$ defines a locally convex topology in $G_M(U)$. 
Note that if $(\varepsilon_m)_{m=1}^{\infty}$ is an arbitrary decreasing to zero sequence of numbers $\varepsilon_m > 0$ then the topology in $G_M(U)$ can be defined also by the system of norms
$$ 
p_m(f)= \sup_{x \in V, \alpha \in {\mathbb Z_+^n}} \frac {\vert (D^{\alpha}f)(x)\vert (1+\Vert x \Vert)^m}{\varepsilon_m^{\vert \alpha \vert} M_{\vert \alpha \vert}} \ , \ m \in {\mathbb N}.
$$ 
Obviously, $G_M(U)$ is continuously embedded in $S(U)$.

We consider a problem of description  of the strong dual space to the spaces $G_M (U)$ in terms of the Fourier-Laplace transform of continuous linear functionals on $G_M (U)$. 

Detailed consideration of this problem depends on additional conditions on the sequence $M$.  
J.W. de Roever [2] studied this problem under the following assumptions on  $M$: 

1). $M_k^2 \le M_{k-1}M_{k+1},  \ \ \forall k \in {\mathbb N}$; 

2). $ \exists H_1 > 1 \ \exists H_2 > 1$  \ $\forall k, m \in {\mathbb Z_+}$ \ \  
$
M_{k+m} \le H_1 H_2^{k+m} M_k M_m;
$

3). $\exists A >0$ \ $\forall m \in {\mathbb N}$ \ 
$
\displaystyle \sum_{k=m+1}^{\infty}\frac {M_{k-1}}{M_k} \le A m \frac {M_m}{M_{m+1}} \ .
$ 

In this case the sequence $M$ is not quasianalytic. Also from conditions 1) and 3) it follows that there exist numbers $h_1, h_2>0$ such that 
$$
M_k \ge h_1 h_2^k k!, \  
\forall k \in {\mathbb Z_+}.
$$

The same problem was considered in [3] under more weak restrictions on $M$. Namely, conditions 2) and 3) were replaced with the following conditions:

$2)'$. $ \exists H_1 > 1 \ \exists H_2 > 1$ \  $\forall k \in {\mathbb Z_+}$  \
$
M_{k+1} \le H_1 H_2^k M_k;
$
      
$3)'$.  $ \exists Q_1 >0 \ \exists Q_2>0$ \  $\forall k \in {\mathbb Z_+}$  \
$
M_k \ge Q_1 Q_2^k k!.
$ 

But the proof of the Paley-Winer type theorem for ultradistributions was given for function  $b$ satisfying the following condition: there is a positive number $r$ such that if $y_1, y_2 \in C$ and $\Vert y_2 - y_1 \Vert\leq 1$
then 
$
\vert b(y_2)-b(y_1)\vert\leq r.
$
 
In this paper we consider the problem in general case. As in [3] we use an approach from B.A. Taylor's paper [4]. 

{\bf 1.2. Definitions and notations}. For $u=(u_1, \ldots , u_m) \in {\mathbb R}^m ({\mathbb C}^m)$, $v=(v_1, \ldots , v_m) \in {\mathbb R}^m ({\mathbb C}^m)$ 
let $<u, v> = u_1v_1 + \cdots + u_m v_m$, $\Vert u \Vert$ be the Euclidean norm in ${\mathbb R}^m ({\mathbb C}^m)$. 

For $z\in\mathbb C^m$, $R>0$ let $B_R (z)$ be a ball in $\mathbb C^m$ of a radius $R$ with the center at the point $z$. Let $\nu_m (R) = \nu_m(1)R^{2m}$ be a volume of $B_R(z)$. 

If $\Omega \subset {\mathbb R}^m$ then $T_{\Omega} = {\mathbb R^m} + i \Omega $. If $\varGamma$ is a cone in $\mathbb R^m$ then 
as usual $pr \varGamma$ is an intersection of $\varGamma$ with the unit sphere.

If $\Omega \subset {\mathbb R}^m$ ($\Omega \subset {\mathbb C}^m $) then the distance from 
$x \in \Omega$ ($z \in \Omega$) to the boundary of $\Omega$ is denoted by $\Delta_{\Omega}(x) (\Delta_{\Omega}(z))$.

For a locally convex space $X$ let $X'$ be the space of linear continuous functionals on $X$ and let $X^*$ be the strong dual space. 

For an open set $\Omega$ in ${\mathbb C}^m$ \ $H(\Omega)$ is a space  of holomorphic functions in $\Omega$, $psh(\Omega)$ is a family of plurisubharmonic  functions in $\Omega$.

If $L=(L_k)_{k=0}^{\infty} $ is a sequence of numbers 
$L_k >0 $ with $L_0=1$ such that 
$$
\displaystyle \lim  \limits_{k \to \infty} \frac {\ln L_k}{k} = + \infty
$$ 
then define a function 
$\omega_L$ on $[0, \infty)$:
$
\omega_L(r) = \displaystyle \sup_{k \in {\mathbb Z_+}} \ln\frac {r^k} {L_k}$ for $r>0$, $\omega_L(0) = 0$. 

{\bf 1.3. Main result}. Throughout this article the non-decreasing sequence $M$ with $M_0=1$ satisfies the following conditions:

$i_1)$. $M_k^2 \le M_{k-1}M_{k+1}, \ k \in {\mathbb N}$; 

$i_2)$. $\exists H_1 > 1 \ \exists H_2 > 1 \ \forall k \in {\mathbb Z_+}$ \ 
$
M_{k+1} \le H_1 H_2^k M_k ;
$
      
$i_3)$. $\exists Q_1 >0 \ \exists  Q_2>0 \ \forall k \in {\mathbb Z_+}$ \ 
$
M_k \ge Q_1 Q_2^k k!.
$ 

Let $S(U)$ be a space of $C^{\infty}(U)$-functions $f$ 
such that for each $p \in {\mathbb Z_+}$
$$
{\Vert f \Vert}_{p, U} = \sup_{x \in V, \vert \alpha \vert \leq p} \vert (D^\alpha f)(x) \vert (1 + {\Vert x \Vert})^p  < \infty. 
$$
Let $S_p(U)$ be a completion of $S(U)$ in the norm ${\Vert \cdot \Vert}_{p, U}$.  
Endow $S(U)$ with a topology of projective limit of spaces $S_p(U)$. Note that $S^*(U)$ is topologically isomorphic to the space of tempered distributions with support in $U$ [5].

If decreasing to zero sequence $(\varepsilon_m)_{m=1}^{\infty}$ of numbers $\varepsilon_m$ is chosen then for brevity denote $\omega_M(\frac{r}{\varepsilon_m})$ by $\omega_m(r), \ r \ge 0$. 

For each $m \in {\mathbb N}$ let
$$
V_{b, m}(T_C) = \{f \in H(T_C): N_m(f) = \sup_{z\in T_C} 
\frac
{\vert f(z)\vert e^{-b(y)}} 
{(1 + \Vert z \Vert)^m(1 + \frac {1} {\Delta_C(y)})^m} < \infty\}, 
$$
$$
H_{b, m}(T_C) = \{ F \in H(T_C): {\Vert F \Vert}_m = \sup_{z \in T_C} 
\frac 
{\vert F(z) \vert}
{e^{b(y) + \omega_m (\Vert z \Vert)} (1 + \frac {1} {\Delta_C (y)})^m} < \infty \} \ , 
$$
where $z = x +i y, x \in {\mathbb R}^n, y \in C$. 

Put $H_{b, M}(T_C) = \bigcup_{m=1}^\infty H_{b, m}(T_C)$, $V_b(T_C) = \bigcup_{m=1}^\infty V_{b, m}(T_C)$.
With usual operations of addition and multiplication by complex numbers $H_{b, M}(T_C)$ and $V_b(T_C)$ are linear spaces. 
Supply $H_{b, M}(T_C)$ ($V_b(T_C)$) with the topology of inductive limit of spaces $H_{b, m}(T_C)$ ($V_{b, m}(T_C)$). 

{\it The Fourier-Laplace transform of a functional} $\varPhi \in S^*(U)$ ($\varPhi \in G_M^*(U)$) is defined by a formula 
$$\hat \varPhi(z) = (\varPhi, e^{i<\xi, z>}), \ z \in T_C.$$ 

In [3] the following theorem was proved. 

{\bf Theorem A.} {\it The Fourier-Laplace transform ${\cal F}: S^*(U) \to V_b(T_C)$  establishes a topological isomorphism between the  spaces $S^*(U)$ and $V_b(T_C)$.}

If $b(y) = a \Vert y \Vert ( a \ge 0)$ then it is a result of V.S. Vladimirov [1].

Let \{$C_{\varepsilon} \}_{\varepsilon > 0}$ be a family of open convex subcones of $C$ such that  
if  ${\varepsilon}_1 < {\varepsilon}_2$ then 
$pr \ \overline {C_{\varepsilon_2}} \subset  pr C_{\varepsilon_1}$ and $\displaystyle \cup_{\varepsilon > 0} 
C_{\varepsilon} = C$. 
Let $H_b(T_C)$ be a projective limit of the spaces $V_{b_{\varepsilon}}(T_C)$, where for $ y  \in {\overline C}$ and $\varepsilon > 0$ \ $b_{\varepsilon}(y) = b(y) + \varepsilon \Vert y \Vert$. Let $H^+_b(T_C)$ be a projective limit of the spaces $H_{b_{\varepsilon}}(T_{C_{\varepsilon}})$. 

We have three useful consequences from theorem A. 

{\bf Corollary 1}. {\it 
The Fourier-Laplace transform ${\cal F}: S^*(U) \to H_b(T_C)$  establishes a topological isomorphism between the  spaces $S^*(U)$ and $H_b(T_C)$.
}

{\bf Corollary 2}. {\it 
The Fourier-Laplace transform ${\cal F}: S^*(U) \to H^+_b(T_C)$  establishes a topological isomorphism between the  spaces $S^*(U)$ and $H^+_b(T_C)$.
}

{\bf Corollary 3}. {\it 
The spaces $V_b(T_C)$, $H_b(T_C)$ and $H^+_b(T_C)$ coincide.
}

The Corollary 2 is known [2].

The aim of the paper is to give  a full proof of the following theorem.
%The following results on description of dual spaces are obtained in the paper.

\begin{theorem}
The Fourier-Laplace transform establishes a topological isomorphism between the  spaces $G_M^*(U)$ and $H_{b, M}(T_C)$.
\end{theorem} 

{\bf Remark 1}. The definition of the spaces $G_M(U)$ and $H_{b, M}(T_C)$ does not depend on a choice of the sequence $(\varepsilon_m)_{m=1}^{\infty}$. So we put ${\varepsilon}_m = H_2^{-m}, \ 
m \in {\mathbb N}$. Thus, $\omega_m(r)= \omega_M(H_2^m r), \ r \ge 0$.

\begin{center}
{\bf 2. Auxiliary results} 
\end{center}

In the proof of theorem 1 the following results will be used.

\begin{lemma} 
For  $y \in C, m \in {\mathbb N}$ let 
$$
g(\xi) = -<\xi, y> + m \ln(1 + \Vert\xi\Vert), \ \xi \in {\mathbb R}^n.
$$ 
Then there exists a number $d>0$ not depending on $y$ such that
$$
\sup \limits_{\xi \in U}g(\xi)
\le 
b(y) + d m + 3m \ln \left(1 + \frac {1}{\Delta_C(y)}\right) + 2m \ln (1 + \Vert y \Vert).
$$
\end{lemma}

Lemma 1 was proved in [3]. 

Due to a special choice of a sequence $(\varepsilon_m)_{m=1}^{\infty}$ the following lemma holds ([6, lemma 2]). 

\begin{lemma} 
For each $ m, p \in {\mathbb N} $ there exists a constant $ c \ge 0 $ such that
$$
w_m(r) + p \ln(1+ r) \le w_{m+p}(r) + c, \ \ r \ge 0 .
$$
\end{lemma}

Also we need a theorem of  J.W. de Roever [2].

{\bf Theorem R.} {\it Let a $n-k$ dimensional hyperplane in ${\mathbb C}^n$ be given by linear functions
$\theta_1 = s_1(\theta_{k+1},
\ldots , \theta_n), \ldots , \theta_k = s_k(\theta_{k+1}, \ldots ,
\theta_n)$ or shortly $w=s(z), w \in {\mathbb C}^k, z \in
{\mathbb C}^{n-k}$.

Let $\Omega_1 \subset \Omega_2 \subset \Omega$ be pseudoconvex domains in ${\mathbb C}^n$ such that an
$\varepsilon > 0$-neighborhood of  $\Omega_1$ with respect to closed polydiscs in the first $k$ coordinates is contained in
$\Omega_2$, i.e.,

$ \{\theta =(\theta_1, \ldots , \theta_n):
\vert \theta_j - \theta_j^{0} \vert \le \varepsilon, j = 1, \ldots
, k; \theta_j = \theta_j^{0}, j = k+1, \ldots , n; \theta^{0} =
(\theta_1^{0}, \ldots ,  \theta_n^{0}) \in \Omega_1 \} \subset
\Omega_2. $

Furthermore, let $\varphi$ be a plurisubharmonic function on $\Omega$ and for
$\theta \in \Omega_1$ let $\varphi_{\varepsilon}(\theta) = \max
\{\varphi (\theta_1 + \xi_1, \ldots , \theta_n + \xi_n: \vert
\xi_j \vert \le \varepsilon, j =1, \ldots , k \}$.

Finally, let
$$
\Omega' = \{z \in {\mathbb C}^{n-k}: (s(z), z) \in \Omega \}, \
\Omega_j' = \{z \in {\mathbb C}^{n-k}: (s(z), z) \in \Omega_j \},
j =1, 2.
$$

And let $\tilde \varphi (z)$ be the function in $\Omega'$ given by $\tilde \varphi (z) =  \varphi(s(z), z), \ z \in  \Omega'$.

Then for a given holomorphic in $\Omega'$ function  $f$ there exists a function $F$ holomorphic in $\Omega_1$ such  that
$F(s(z), z)=f(z), \ z \in  \Omega'$,  and for some $K>0$ depending only on $k$ and $s_1, \ldots , s_k$,
$$
\int_{\Omega_1} \frac
{F(\theta)\exp(-\varphi_{\varepsilon}(\theta))} {(1 + \vert \theta
\Vert^2)^{3k}} \ d  \lambda_n (\theta) \le K \varepsilon^{-2k}
\int_{\Omega_2'} \vert f(z) \vert^2 e^{-\tilde \varphi (z)} \ d
\lambda_{n-k} (\theta),
$$
where $\lambda_n$ and $\lambda_{n-k}$ denote the Lebesgue measure in ${\mathbb
C}^n$ or ${\mathbb C}^{n-k}$, respectively, if $f$ is such that the right hand side is finite. $F$ depends besides on $f$ also on $\Omega_1, \varepsilon, \varphi$.
}

\begin{lemma} 
Let $\varGamma$ be an open convex cone in ${\mathbb R}^n$ with the apex at the origin.
Let $h$ be a convex continuous positively homogeneous function of degree 1 on ${\overline \varGamma}$. Then for each $\varepsilon > 0$ there exists a constant $A_{\varepsilon} > 0$ such that for $y_1, y_2 \in \varGamma$ satisfying the inequality $\Vert y_2 - y_1 \Vert \le 1$ 
$$
\vert h(y_2) - h(y_1) \vert \le \varepsilon \Vert y_1 \Vert +  \varepsilon \Vert y_2 \Vert + A_{\varepsilon} .
$$
\end{lemma} 

The proof is easy and so it is omitted. 

Let for brevity for $z \in T_C$  $d(z)= \Delta_{T_C}(z)$.

\begin{theorem} 
Let  function $S \in H({\mathbb C^n \times T_C})$ for some $m \in {\mathbb N}$  on ${\mathbb C^n \times T_C}$ satisfies the inequality 
$$
\vert S(z, \zeta) \vert \leq \exp (\omega_m(\Vert z\Vert) + b(Im \zeta) + m\ln (1 + \Vert\zeta\Vert) + m \ln (1 + \frac{1}{d(\zeta)})),
$$
and $S(\zeta, \zeta) = 0, \ \ \zeta\in T_C$.

Then there exist functions $S_1, \ldots , S_n \in H({\mathbb C^n \times T_C})$ and a number $p \in {\mathbb N}$ such that:

a) $ S(z, \zeta) = \displaystyle\sum_{j=1}^n S_j (z, \zeta)(z_j - \zeta_j),\ \ (z, \zeta)\in {\mathbb C^n \times T_C}$;
 
b) for each $\varepsilon > 0$ there exists a number $a_{\varepsilon} > 0 $ such that for $j=1, \ldots , n$ and 
$(z, \zeta) \in {\mathbb C^n \times T_C}$
$$
\vert S_j(z, \zeta) \vert \leq a_{\varepsilon} \exp ( \omega_p(\Vert z\Vert) + b(Im \zeta)+ \varepsilon \Vert Im \zeta \Vert  + p \ln (1 + \frac{1} {d(\zeta)}) + p\ln (1 + \Vert\zeta\Vert)).
$$
\end{theorem} 

{\bf Proof}. Let for $\zeta \in T_C, j\in {\mathbb N}$ 
$$
\varphi_j (\zeta) = b(Im \zeta) + j \ln((1+\frac{1}{d(\zeta)})(1+\Vert \zeta \Vert )),
$$
for $\zeta \in T_C, z_1, \ldots,z_k \in {\mathbb C}, l, k, j \in {\mathbb N}$
$$
h_{l,k,j} (z_1,\ldots,z_k,\zeta) = \omega_l (\Vert (z_1+\zeta_1,\ldots,z_k+\zeta_k,\zeta_{k+1},\ldots,\zeta_n)\Vert) + \varphi_j (\zeta).
$$
Functions $h_{l,k,j}$ are plurisubharmonic in ${\mathbb C^k} \times T_C$. 

Put $L(z, \zeta) = S(z+\zeta, \zeta)$, \ $z \in {\mathbb C^n}, \zeta\in T_C.$ 
By assumption, 
$$
\vert L(z, \zeta)\vert \leq \exp(h_{m,n,m} (z, \zeta)), \  (z, \zeta)\in {\mathbb C^n \times T_C}, \eqno (1)
$$ 
and 
$L(0, \zeta) = S(\zeta, \zeta) = 0$ for $\zeta \in T_C$.

We now pass to construction of functions $S_1, \ldots , S_n$.
Let 
$L_1 (z_1,\zeta) = L(z_1,0,\ldots,0,\zeta)$,  \ $z_1\in {\mathbb C}$, $\zeta\in T_C ;$
$L_2 (z_1,z_2,\zeta) = L(z_1,z_2,0,\ldots,0,\zeta)$,\ $z_1,z_2 \in {\mathbb C}$, $ \zeta\in T_C;  \ldots $ , 
$L_n(z, \zeta) = L(z,\zeta), \ z\in {\mathbb C^n}, \zeta\in T_C.$
In view of (1) we have 
$$
\vert L_k(z,\zeta)\vert \leq \exp(h_{m,k,m}(z,\zeta)), \ z \in {\mathbb C^k}, \zeta \in T_C. \eqno (2)
$$
Since $L_1(0,\zeta)=0 \ \ \forall\zeta\in T_C$, then the function 
$$
\psi_1^{(1)}(z_1,\zeta) = \frac{L_1(z_1,\zeta)}{z_1}
$$
is holomorphic in ${\mathbb C}\times T_C$. Let us estimate the growth of the function $\psi_1^{(1)}$. 
If
$\vert z_1 \vert \geq 1$, then $\forall\zeta\in T_C$
$$\vert\psi_1^{(1)}(z_1, \zeta)\vert \leq  \vert L_1(z_1, \zeta)\vert \leq \exp(h_{m,1,m}(z_1, \zeta)),
$$
If $\vert z_1\vert <1$, then $\forall\zeta\in T_C$
$$
\vert\psi_1^{(1)}(z_1, \zeta)\vert \leq \max_{\vert t_1\vert =1}\vert \psi_1^{(1)} (t_1, \zeta)\vert = \max_{\vert t_1\vert =1}\vert L_1(t_1, \zeta)\vert. 
$$
In view of the inequality (2) we have
$$
\vert \psi_1^{(1)}(z_1,\zeta)\vert \leq  \exp ({\displaystyle\max_{\vert t_1\vert =1} h_{m, 1, m}(t_1, \zeta)}).
$$
Recall that for each $m \in {\mathbb N}$ there exists a number $K_m > 0$ (see [6], lemma 1 (in the proof of this lemma only the conditions $i_1), i_3)$ are used) such that
$$
\vert \omega_m (r_2) - \omega_m(r_1)\vert \leq K_m\vert r_2 - r_1\vert,\ \ r_1, r_2 \geq 0.  \eqno (3)
$$
So, using (3) we obtain
$$
\vert \psi_1^{(1)}(z_1,\zeta)\vert\leq A_1 e^{h_{m,1,m}(z_1, \zeta)}, \ z_1 \in {\mathbb C}, \zeta \in T_C,
$$
where $A_1 = e^{2K_m}$. 
In addition $L_1(z_1,\zeta) = \psi_1^{(1)}(z_1, \zeta)z_1$  everywhere in ${\mathbb C}\times T_C$. 

Let for $z_1\in~{\mathbb C}$, $\zeta \in T_C$
$$
\varphi^{(1)}(z_1, \zeta)=2 h_{m, 1, m}(z_1, \zeta) + (2n+3) \ln(1+\Vert(z_1,\zeta)\Vert).
$$
Then
$$
\int \limits_{{\mathbb C}\times T_C} \vert \psi_1^{(1)}(z_1,\zeta)\vert^2 e^{-\varphi^{(1)}(z_1, \zeta)} \ d\lambda_{n+1}(z_1, \zeta) < \infty .
$$
Suppose that for $k=2, \ldots, n$ there are found functions $\psi_j^{(k-1)}$ holomorphic in ${\mathbb C}^{k-1} \times T_C$ ($j=1, \ldots, k-1$)
such that
$$
L_{k-1}(z_1, \ldots, z_{k-1}, \zeta) = \displaystyle \sum_{j=1}^{k-1} \psi_j^{(k-1)}(z_1,\ldots,z_{k-1}, \zeta)z_j
$$
and for some $d_{k-1}\in {\mathbb N}$  and for all $j= 1, \ldots , k-1$
$$
\int \limits_{{\mathbb C}^{k-1} \times T_C} 
\vert \psi_j^{(k-1)}(z, \zeta)\vert^2 e^{-(2 h_{m, k-1,m}(z, \zeta) + d_{k-1} \ln(1+\Vert(z, \zeta)\Vert))}  < \infty . \eqno (4)
$$
Now we use theorem R and notations of this theorem. Let us consider in ${\mathbb C^k}\times {\mathbb C^n}$ the subspace $$
W_k=\{(z_1,z_2,\ldots,z_k,\zeta):z_k=0,z_1,\ldots,z_{k-1} \in {\mathbb C}, \zeta\in {\mathbb C^n}\}
$$
of codimension 1. 
Put $\Omega=\Omega_1=\Omega_2={\mathbb C^k}\times T_C$. 
Then 
$$
\Omega^{'}=\{(z_1,\ldots, z_{k-1},\zeta):(z_1, \ldots, z_{k-1},0,\zeta)\in {\mathbb C^k}\times T_C\}
={\mathbb C^{k-1}}\times T_C.
$$ 
So $\Omega_1^{'}=\Omega_2^{'}=\Omega^{'}$.
Put 
$$
\varphi^{(k)}(z, \zeta)=2 h_{m, k, m}(z,\zeta) + d_{k-1}\ln(1+\Vert(z,\zeta)\Vert),
$$
where 
$z\in {\mathbb C^k}, \zeta\in T_C.
$
In this theorem for $\varepsilon$ we take 1. 
Let
$$
\varphi_1^{(k)}(z_1, \ldots, z_k, \zeta) = \displaystyle \max_{\vert \xi_k \vert \leq1} \varphi^{(k)}(z_1, \ldots, z_{k-1}, z_k+\xi_k, \zeta).
$$ 
Then
$$
\widetilde{\varphi}^{(k)}(z_1, \ldots, z_{k-1}, \zeta)= \varphi^{(k)}(z_1, \ldots, z_{k-1}, 0, \zeta)=
$$
$$
=2 h_{m, k-1, m} (z_1, \ldots, z_{k-1}, \zeta) + d_{k-1} \ln(1+\Vert(z_1, \ldots, z_{k-1}, \zeta)\Vert)],  
$$
where $z_1,\ldots,z_{k-1} \in {\mathbb C}, \zeta\in T_C$. 
In view of (4) 
$$
\int_{{\mathbb C^{k-1}}\times T_C}{\vert \psi_j^{(k-1)} (z, \zeta)\vert}^2 e^{-\widetilde{\varphi}^{(k)}(z, \zeta)} \ 
d\lambda_{k-1+n}(z, \zeta) <\infty, \ j=1, \ldots, k-1.
$$
By Roever's theorem  there exists a function $\psi_j^{(k)}$ ($j=1,\ldots,k-1 $) holomorphic in ${\mathbb C^k} \times T_C$ such that
$$
\psi_j^{(k)}(z_1, \ldots, z_{k-1}, 0, \zeta) = \psi_j^{(k-1)}(z_1, \ldots, z_{k-1}, \zeta), \ z_1, \ldots, z_{k-1} \in {\mathbb C}, \zeta\in T_C, $$
and for some $K_j^{(k-1)}>0$ not depending on $\psi_j^{(k-1)}$
$$
\int_{{\mathbb C^k} \times T_C} \frac {{\vert\psi_j^{(k)}(z, \zeta)\vert}^2 e^{-\varphi_1^{(k)}(z, \zeta)}} {(1 + {\Vert(z, \zeta)\Vert}^2)^3} \ d \lambda_{n+k}(z, \zeta)\leq
$$
$$
\leq 
K_j^{(k-1)}\int_{{\mathbb C^{k-1}}\times T_C}{\vert\psi_j^{(k-1)}(z', \zeta)\vert}^2e^{-\widetilde{\varphi}^{(k)}
(z', \zeta)}\ d\lambda_{n+k-1}(z', \zeta).
$$
Using (3) we have
$$
\int_{{\mathbb C^k} \times T_C} \frac {{\vert\psi_j^{(k)}(z, \zeta)\vert}^2 e^{-\varphi^{(k)}(z, \zeta)}} {(1 + {\Vert(z, \zeta)\Vert}^2)^3} \ d \lambda_{n+k}(z, \zeta) < \infty \ , j=1,\ldots,k-1.
$$
Now consider the function 
$$
V_k(z, \zeta)=L_k(z, \zeta) -\psi_1^{(k)}(z, \zeta) z_1 - \ldots - \psi_{k-1}^{(k)}(z, \zeta)z_{k-1}, \ z \in\mathbb C^k,  \zeta \in T_C.
$$
It is holomorphic in ${\mathbb C^k} \times T_C$.
Note that
$$
\int_{{\mathbb C^k} \times T_C} \frac {{\vert V_k(z, \zeta) \vert}^2 e^{-\varphi^{(k)}(z, \zeta)}} {(1 + {\Vert(z, \zeta)\Vert}^2)^4} \ d \lambda_{n+k}(z, \zeta) < \infty .
$$
Since for $z_1, \ldots , z_{k-1} \in \mathbb C, \zeta \in  T_C$ \ 
$
V_k(z_1, \ldots,z_{k-1}, 0, \zeta)= 0,
$
then the function
$$
\psi_k^{(k)}(z, \zeta)=\frac{V_k(z, \zeta)} {z_k} 
$$
is holomorphic in ${\mathbb C^k} \times T_C$.

Put $\gamma^{(k)}(z, \zeta) = \varphi^{(k)}(z, \zeta) + 8 \ln(1+\Vert(z,\zeta)\Vert)$.
Obviously,
$$
\int_{{\mathbb C}^{k-1} \times \{\vert z_k \vert > 1 \} \times T_C} 
\vert \psi_k^{(k)}(z, \zeta) \vert^2 e^{-\gamma^{(k)}(z, \zeta)} \ d \lambda_{n+k}(z, \zeta) < \infty .
$$
Let $\vert z_k \vert \le 1$. Then for $z =(z', z_k) \in {\mathbb C}^k, \zeta \in T_C$
$$
\vert \psi_k^{(k)}(z', z_k, \zeta)\vert^2 \le 
\max \limits_{\vert \eta \vert =1} \vert V_k(z', \eta, \zeta)\vert^2 \le 
\frac {1}{\pi} \int_{\vert t \vert \le 2}\vert V_k(z', t, \zeta)\vert^2 \ d \lambda_1(t) .
$$
Hence,
$$
\vert \psi_k^{(k)}(z', z_k, \zeta)\vert^2 e^{-\gamma^{(k)}(z', z_k,  \zeta)} \le 
\frac {1}{\pi} \int_{\vert t \vert \le 2}\vert V_k(z', t, \zeta)\vert^2 e^{-\gamma^{(k)}(z', z_k,  \zeta)} \ d \lambda_1(t) .
$$
Further,
$$
\int_{{\mathbb C^{k-1}} \times T_C}
\vert \psi_k^{(k)}(z', z_k, \zeta)\vert^2 e^{-\gamma^{(k)}(z', z_k,  \zeta)} \ d \lambda_{n+k-1}(z', \zeta) \le 
$$
$$
\le 
\frac {1}{\pi} \int_{\vert t \vert \le 2} \int_{{\mathbb C^{k-1}} \times T_C} \vert V_k(z', t, \zeta)\vert^2 e^{-\gamma^{(k)}(z', z_k,  \zeta)} \ d \lambda_{n+k-1}(z', \zeta)\ d \lambda_1(t) .
$$
Since for $z =(z', z_k) \in {\mathbb C}^k, \zeta \in T_C$, $\vert z_k \vert \le 1, \vert t \vert \le 2$
$$
\vert \gamma^{(k)}(z', t,  \zeta) - \gamma^{(k)}(z', z_k,  \zeta) \vert \le 6 K_{m_{k-1}} + 3 (d_{k-1} + 8),
$$
then, putting 
$b_{k-1} = 6 K_m + 3 (d_{k-1} + 8)$, 
we have
$$
\int_{{\mathbb C^{k-1}} \times T_C}
\vert \psi_k^{(k)}(z', z_k, \zeta)\vert^2 e^{-\gamma^{(k)}(z', z_k,  \zeta)} \ d \lambda_{n+k-1}(z', \zeta)
\le
$$
$$
\le 
\frac {e^{b_{k-1}}}{\pi} \int_{\vert t \vert \le 2} \int_{{\mathbb C^{k-1}} \times T_C} \vert V_k(z', t, \zeta)\vert^2 e^{-\gamma^{(k)}(z', t,  \zeta)} \ d \lambda_{n+k-1}(z', \zeta)\ d \lambda_1(t) .
$$
From here it follows that
$$
\int_{\vert z_k \vert \le 1} \int_{\mathbb C^{k-1} \times T_C}
\vert \psi_k^{(k)}(z', z_k, \zeta)\vert^2 e^{-\gamma^{(k)}(z', z_k,  \zeta)} \ d \lambda_{n+k-1}(z', \zeta) 
\ d \lambda_1(z_k)
\le
$$
$$
\le 
e^{b_{k-1}} \int_{\vert t \vert \le 2} \int_{{\mathbb C^{k-1}} \times T_C} \vert V_k(z', t, \zeta)\vert^2 e^{-\gamma^{(k)}(z', t,  \zeta)} \ d \lambda_{n+k-1}(z', \zeta)\ d \lambda_1(t) .
$$
Finally, 
$$
\int_{{\mathbb C}^k \times T_C}
\vert \psi_k^{(k)}(z, \zeta)\vert^2 e^{-\gamma^{(k)}(z, \zeta)} \ d \lambda_{n+k}(z, \zeta) 
\le
$$
$$
\le 
2 e^{b_{k-1}} \int_{{\mathbb C}^k \times T_C} \vert V_k(z, \zeta)\vert^2 e^{-\gamma^{(k)}(z, \zeta)} \ d \lambda_{n+k}(z, \zeta).
$$
Thus, there are functions $\psi_j^{(k)} \in H({\mathbb C^k} \times T_C)$ ($j=1, \ldots, k$) such that
$$
L_k(z_1, \ldots, z_k, \zeta) = \displaystyle \sum_{j=1}^k \psi_j^k(z_1, \ldots, z_k, \zeta)z_j .
$$
and
$$
\int_{{\mathbb C^k} \times T_C} 
\vert \psi_j^{(k)}(z, \zeta) \vert^2 e^{-\gamma^{(k)}(z, \zeta)} \ d \lambda_{n+k}(z, \zeta) < \infty ,
$$
Putting $d_k = d_{k-1} + 8$, we have
$$
\int_{{\mathbb C^k} \times T_C} 
\vert \psi_j^{(k)}(z, \zeta) \vert^2 e^{-(2 h_{m, k, m}(z,\zeta) + d_k \ln(1+\Vert(z,\zeta)\Vert))} \ d \lambda_{n+k}(z, \zeta) < \infty , \ j=1, \ldots, k.
$$
If $k=n$, then we have
$$
L_n(z, \zeta)=L(z, \zeta)=\psi_1^{(n)}(z, \zeta)z_1+\ldots+\psi_n^{(n)}(z, \zeta)z_n
$$
and for  some $ d_n \in {\mathbb N}$ for all $j=1, \ldots, n$
$$
C_{j, n}:=\int_{{\mathbb C^n} \times T_C} 
\vert \psi_j^{(n)}(z, \zeta) \vert^2 e^{-(2 h_{m, n, m}(z,\zeta) + d_n \ln(1+\Vert(z,\zeta)\Vert))} \ d \lambda_{n+k}(z, \zeta) < \infty . 
$$

Let us obtain uniform estimates of functions $\psi_j^{(n)}$. Let $(z, \zeta) \in {\mathbb C^n \times T_C}$. Put $R=\min (1, \frac {d(\zeta)}{4})$.  
Let $(t, u) \in {\mathbb C^n \times T_C}$ belongs to the ball $B_R((z, \zeta))$. Then 
$$
\vert h_{m, n, m}(t,u)- h_{m, n, m}(z,\zeta)\vert \le \vert w_m (\Vert (t+u)\Vert) - w_m (\Vert (z+\zeta)\Vert) \vert + 
\vert \varphi_m(u) - \varphi_m(\xi) \vert .
$$
The first term from the right is estimated as follows (using (3)):
$$
\vert w_m (\Vert (t+u)\Vert) - w_m (\Vert (z+\zeta)\Vert) \vert \le K_m (\Vert t - z \Vert + \Vert u - \zeta \Vert) \le 2K_m. 
$$
Furthermore,
$$
\vert \varphi_m(u) - \varphi_m(\zeta) \vert  \le \vert b(Im u) - b(Im \zeta) \vert + 
m \vert \ln(1+\frac{1}{d(u)}) -  \ln(1+\frac{1}{d(\zeta)}) \vert + 
$$
$$
+
m \vert \ln(1+\Vert u \Vert ) - \ln(1+\Vert \zeta \Vert )\vert .
$$
Note that
$$
\left \vert \ln \left(1+\frac{1}{d(u)}\right) -  \ln \left (1+\frac{1}{d(\zeta)}\right) \right \vert = \left \vert \ln 
\frac 
{(1+d(u))d(\zeta)}{(1+d(\zeta))d(u)} \right \vert  \le 
$$
$$
\le 
\left \vert \ln 
\frac 
{4(1+d(u))}{3(1+d(\zeta))} \right \vert  \le \left \vert \ln \left (1 + \frac {d(u) - d (\zeta)}{1 + d (\zeta)} \right) \right \vert + 1 \le 
$$
$$
\le 
\left \vert \frac {d(u) - d (\zeta)}{1 + d (\zeta)}  \right \vert + 1 \le \frac 5 4.
$$
Besides that, 
$$
\vert \ln (1 + \Vert u \Vert) - \ln (1 + \Vert \zeta \Vert) \vert \le \Vert u - \zeta \Vert \le 1.
$$
By lemma 3  for each $\varepsilon > 0$ there is a constant $A_{\varepsilon} > 0$ (not depending on $\zeta$) such that  
$$
\vert b(Im u) - b(Im \zeta) \vert 
\le \varepsilon \Vert Im u \Vert +  \varepsilon \Vert Im \zeta \Vert + A_{\varepsilon}.
$$
So we have
$$
\vert b(Im u) - b(Im \zeta) \vert \le 2 \varepsilon \Vert Im \zeta \Vert + A_{\varepsilon} + \varepsilon .
$$
Thus, for each $\varepsilon > 0$ and for
all $u, \zeta \in T_C$ such that $\Vert u - \zeta \Vert \le R$ 
$$
\vert \varphi_m(u) - \varphi_m(\zeta) \vert  \le 2 \varepsilon \Vert Im \zeta \Vert + A_{\varepsilon} + \varepsilon  + \frac {9m} {4} . 
$$
Let $c_{\varepsilon} = A_{\varepsilon} + \varepsilon  + \frac {9m} {4} + 2K_m$. 
Then for each $(t, u) \in B_R((z, \zeta))$ 
$$
\vert h_{m, n, m}(t,u)- h_{m, n, m}(z,\zeta)\vert \le 
2 \varepsilon \Vert Im \zeta \Vert + c_{\varepsilon} . \eqno (5)
$$
For a plurisubharmonic in $\mathbb C^n \times T_C$ function $ \vert \psi_j^{(n)}(t, u) \vert^2 $ we have
$$
\vert \psi_j^{(n)}(z, \zeta) \vert^2 \leq 
\frac 
{1}{\nu_{2n}(R)}
\displaystyle \int \limits_{B_R((z, \zeta))}
\vert \psi_j^{(n)}(t, u)\vert^2 \ d\lambda_{2n}(t, u) \leq  
$$
$$
\leq 
\frac 
{1}{\nu_{2n}(R)}
\int \limits_{\mathbb C^n\times T_C} \vert\psi_j^{(n)}(t, u)\vert^2 e^{-(2 h_{m, n, m}(t, u) + d_n \ln(1+\Vert(t, u)\Vert))}  
 \ d\lambda_{2n}(t, u) \times
$$
$$
\times e^{\sup \limits_{(t, u) \in B_R((z,\zeta))} (2 h_{m, n, m}(t, u) + d_n \ln(1+\Vert(t, u)\Vert))}  \le 
$$
$$
\le 
\frac 
{C_{j, n}}{\nu_{2n}(R)}
 (2 + \Vert (z,\zeta) \Vert)^{d_n} e^{\sup \limits_{(t, u) \in B_R((z,\zeta))} (2 h_{m, n, m}(t, u))}
$$
Using (5) for each $\varepsilon > 0$ we find a constant $b_{\varepsilon} > 0$ such that
$$
\vert \psi_j^{(n)}(z, \zeta) \vert^2 \leq 
\frac 
{1}{\nu_{2n}(R)}
C_{j, n} b_{\varepsilon} (2 + \Vert (z,\zeta) \Vert)^{d_n}  e^{2 h_{m, n, m}(z,\zeta) + 2 \varepsilon \Vert Im \zeta \Vert} .
$$
Thus, for each $\varepsilon > 0$ there is a constant $M_{\varepsilon} > 0$ such that
$$
\vert \psi_j^{(n)}(z, \zeta) \vert \leq M_{\varepsilon} \left(1 + \frac {1}{d(\zeta)}\right)^{2n}(1 + \Vert (z,\zeta) \Vert)^{d_n}e^{h_{m, n, m}(z,\zeta) + \varepsilon \Vert Im \zeta \Vert}.
$$
From here 
it  follows that for each $\varepsilon > 0$ 
we can find a constant $L_{\varepsilon} > 0$ such that
$$
\vert \psi_j^{(n)}(z, \zeta) \vert \leq L_{\varepsilon} \left(1 + \frac {1}{d(\zeta)}\right)^{2n}
(1 + \Vert (z +\zeta,\zeta)\Vert)^{d_n}(1 + \Vert \zeta \Vert)^{d_n}e^{h_{m, n, m}(z,\zeta) + \varepsilon \Vert Im \zeta \Vert}.
$$
Let $t_n= \max (d_n, 2n), p = m+t_n$. 
Then, by lemma 2 for each $\varepsilon > 0$ we find a constant $a_{\varepsilon} > 0$ such that
$$
\vert \psi_j^{(n)}(z, \zeta) \vert \leq a_{\varepsilon}
e^{h_{p, n, p}(z,\zeta) + \varepsilon \Vert Im \zeta \Vert}.
$$
Put $S_j(z, \zeta)=\psi_j^{(n)}(u-\zeta,\zeta)$. Then $S_j \in H({\mathbb C^n} \times T_C)$, 
$$
S(z, \zeta)=\sum_{j=1}^n S_j(z,\zeta)(z_j-\zeta_j), \ z \in \mathbb C^n, \zeta \in T_C,
$$
and for each $\varepsilon > 0$ and for all $(z, \zeta) \in {\mathbb C^n \times T_C}$
$$
\vert S_j(z,\zeta)\vert\leq a_{\varepsilon}\exp(\omega_p(\Vert z \Vert)+ b(Im \zeta) + \varepsilon \Vert Im \zeta \Vert + p \ln(1+\frac{1}{d(\zeta)}) + p \ln(1+\Vert \zeta \Vert)).
$$

\begin{center}
{\bf 3. Proof of theorem 1} 
\end{center}

Let
$$
C_m(U) = \{f\in C(U): \widetilde{p}_m (f) = 
\sup_{x \in U} \vert f(x) \vert (1+\Vert x \Vert)^m < \infty \}, \  m \in {\mathbb N}.
$$
By standard scheme ([4, Propositions 2.10, 2.11, Corollary 2.12]) one can prove the following lemma. 

\begin{lemma}
Let $T\in G_M'(U)$ and numbers $c>0$, $m\in\mathbb N$ be such that
$$
\vert (T, f ) \vert \leq c p_m(f), \ f \in G_M(U).
$$
 
Then there exist functionals $T_{\alpha}\in C_m'(U)$ ($\alpha\in\mathbb Z_+^n$) such that
$$
\vert (T_{\alpha},f) \vert\leq\frac{C\widetilde{p}_m(f)}  {\varepsilon_m^{\vert\alpha\vert}M_{\vert\alpha\vert}} \ , \ f \in C_m(U),
$$ 
and 
$
(T, f)=\displaystyle\sum_{\vert\alpha\vert\geq 0}(T_{\alpha}, D^{\alpha}f ), \ 
f \in G_M(U). 
$
\end{lemma}

\begin{lemma} 
Let $S \in G'_M(U)$. Then $\hat S \in H_{b, M}(T_C)$. 
\end{lemma}

{\bf Proof}. First note that for each $z=x + i y \in T_C$ ($x \in {\mathbb R}^n, y \in C$) 
the function $f_z(\xi) = e^{i<\xi ,z>}$ belongs to  $G_M(U)$. 
Indeed, for each $m \in {\mathbb N}$ 
$$
p_m(f_z) = 
\sup_{\xi \in V, \alpha \in {\mathbb Z_+^n}} \frac {\vert (iz)^\alpha e^{i<\xi, z>} \vert(1 + \Vert \xi \Vert)^m} {\varepsilon_m^{\vert \alpha \vert} M_{\vert \alpha \vert}} \le
$$
$$
\le 
\sup_{\alpha \in {\mathbb Z_+^n}} 
\frac {{\Vert z \Vert}^{\vert \alpha \vert}}  {\varepsilon_m^{\vert \alpha \vert} M_{\vert \alpha \vert} }
\sup \limits_{\xi \in V}
e^{ -<\xi, y> + m \ln(1 + \Vert \xi \Vert)} =
e^{\omega_m(\Vert z \Vert) +\sup \limits_{\xi \in V}( -<\xi, y> + m \ln(1 + \Vert \xi \Vert))}. 
$$
Using lemma 1 and 2, we get
$$
p_m(f_z) \le A e^{b(y) + \omega_{3m} (\Vert z \Vert)} \left(1 + \frac {1}{\Delta_C(y)}\right)^{3m} \ , \eqno (6)
$$
where $A = e^{d m}$ doesn't depend on $z \in T_C$. 

Now let $S \in G'_M(U)$. It is clear that $\hat S$ is correctly defined on $T_C$. 
Using lemmas 1 and 4, condition $i_3)$ it is easy to see that $\hat S \in H(T_C)$. 

Since there exist numbers $m \in {\mathbb N}$ and $c>0$ such that
$$
\vert (S, f) \vert \leq c p_m(f), \ \ f\in G_M(U),
$$
then (using (6)) we obtain
$$
\vert \hat S(z) \vert \leq  cA e^{b(y) + \omega_{3m} (\Vert z \Vert)} \left(1 + \frac {1}{\Delta_C(y)}\right)^{3m} \ .  
$$
Thus, $\hat S \in H_{b, M}(T_C)$ and lemma is proved.

Obviously, for each $m \in {\mathbb N}$ the embeddings $j_m : H_{b, m}(T_C) \rightarrow H_{b, m+1}(T_C)$ are completely continuous. 
Consequently, $H_{b, M}(T_C)$ is an $(LN^*)$-space. 

{\bf Proof of theorem 1}. By lemma 5 the linear mapping $L: S \in G_M^*(U) \rightarrow \hat S$ acts from $G_M^*(U)$ into $H_{b,M}(T_C)$.

Before showing the coninuity of $L$, note that the topology of the space $G_M^*(U)$ can be described as follows.
Let for each $k \in {\mathbb N}$
$$
W_k = \{ f \in G_M(U):  \ p_k(f) \leq 1 \},  
$$ 
$$
W_k^0 = \{ F \in G_M'(U): \vert (F, f) \vert \leq 1,  \ \forall f \in W_k \}.
$$
Let 
$E_k = \displaystyle \bigcup_{\alpha >0} (\alpha W_k^0) $.
Endow $E_k$ with the norm
$$
q_k(F) = \displaystyle \sup_{f \in W_k} \vert (F, f) \vert,  \ F\in E_k.
$$
Note that
$G_M'(U) = \bigcup_{k=1}^\infty E_k.$ 
Define in $G_M'(U)$ the topology $\lambda$ of the inductive limit of spaces $E_k$. 
Since $G_M(U)$ is an $(M^*)$-space then $G_M(U)$ is Montel space (and hence reflexive). 
But then the strong topology in $G_M'(U)$ coincides with the topology $\lambda$ [7, p. 699-700]. 

Now let $S \in E_m, \ m \in{\mathbb N}$. 
Then 
$$
\vert (S, f) \vert \leq q_m(S), \ \forall f \in W_m. 
$$
Hence,
$$
\vert (S, f) \vert \leq q_m(S) p_m(f), \ f \in G_M(U).
$$ 
From here and (6) we obtain 
$$
\vert \hat S(z) \vert  \leq q_m(S) A e^{b(y) + 
\omega_{3m} (\Vert z \Vert)} \left(1 + \frac {1}{\Delta_C(y)}\right)^{3m} \ ,  
$$\\
where a constant $A>0$ doesn't depend on $z \in T_C$.
Hence, 
$$
{\Vert \hat S\Vert}_{3m} \leq  A q_m(S),  \ S \in E_m  \ (m =1, 2, \ldots),
$$ 
and, consequently, $L$ is continuous.

Now let us prove that $L$ is bijective.

We first show that $L$ is surjective.
Let $F \in H_{b, M}(T_C)$. That is, $F \in H(T_C)$ and for some $c>0, m \in {\mathbb N}$ 
$$
\vert F(z) \vert \leq c e^{b(y) + \omega_m (\Vert z \Vert)} \left (1 + \frac {1} {\Delta_C (y)}\right)^m, \ z \in T_C.  \eqno (7)
$$
Note that for $z = x + iy \in T_C$ \ $d(z) = \Delta_C (y)$. 
Hence, from (7) we have
$$
\int_{T_C}{\vert F(z)\vert}^2 e^{-2(b(y) + \omega_m (\Vert z \Vert) + m\ln (1 + \frac{1}{d(z)}) + (n+1)\ln (1 + \Vert z\Vert^2))} \ 
d \lambda_n (z) <\infty. \eqno (8)
 $$
Let $K = {\mathbb R^n} \times C$. 
Put in theorem R (with the replacement of  $n$ with $2n$) \ $\Omega = \Omega_1 = \Omega_2 = {\mathbb R^{2n}} + iK$. 
Note that $\Omega = \Omega_1 = \Omega_2 = {\mathbb C^n} \times T_C$. 
Since $K$ is a convex domain in 
${\mathbb R^{2n}}$ then ${\mathbb R^{2n}} + iK$ is a domain of holomorphy. 
As a linear subspace in this theorem consider 
$$
W = \{ (z, \xi) \in {\mathbb C^{2n}}: \ \ z_1=\xi_1,\ldots,z_n=\xi_n\}
$$
in ${\mathbb C^{2n}}$ of complex dimension $n$. 
Then
$$
\Omega^{'} = \Omega_1^{'} = \Omega_2^{'} = \{z\in {\mathbb C^n}: (z, z) \in\Omega = {\mathbb C^n}\times T_C\} = T_C.
$$ 
In theorem R as $\varepsilon$ we take 1 and as $\varphi$ we take the function
$$
\varphi(z, \xi) = 2(b(Im \xi)  + \omega_m (\Vert z \Vert) + m \ln (1 + \frac{1}{d(\xi)}) + (n+1)\ln(1 + {\Vert (z, \xi)\Vert}^2)), 
$$
where $z = x +iy \in {\mathbb C^n}, \xi \in T_C$. 
Note that $\varphi \in psh ({\mathbb C^n} \times T_C$) and
$$
\widetilde{\varphi}(z) = 2 (b(y) + \omega_m (\Vert z \Vert) + m \ln (1 + \frac{1}{d(z)}) + (n+1)\ln(1 + 2{\Vert z\Vert}^2)), \ 
z\in T_C.
$$
In view of (8)
$$
\int_{T_C}{\vert F(z)\vert}^2 e^{-\widetilde{\varphi}(z)} \ d \lambda_n(z) < \infty.  
$$
By Roever's theorem there exists a function $\Phi \in H({\mathbb C^n} \times T_C)$ such that $\Phi(z, z) = F(z)$ for $z\in T_C$ and for some $B>0$
$$
\int_{{\mathbb C^n} \times T_C} \frac {{\vert \Phi(z, \xi)\vert}^2 e^{-\varphi_1 (z, \xi)}} {(1 + {\Vert (z, \xi)\Vert}^2)^{3n}} \  d\lambda_{2n} (z, \xi) \leq B \int_{T_C} {\vert F(z)\vert}^2 e^{-\widetilde{\varphi}(z)} \ d \lambda_n (z).
$$
Here 
$$
\varphi_1 (z, \xi) = \displaystyle \max_{\vert t_1 \vert \leq 1, \ldots,\vert t_n\vert\leq 1} \varphi (z_1+t_1, \ldots, z_n+t_n, \xi_1, \ldots, \xi_n).
$$
Using (3) and an inequality
$$
\vert \ln (1 + x_2^2) - \ln (1 + x_1^2) \vert \leq \vert x_2 - x_1 \vert, \ \ x_1,x_2 \in {\mathbb R},
$$
we obtain in ${\mathbb C^n} \times T_C$ \ 
$
\vert \varphi_1 (z, \xi) - \varphi (z, \xi)\vert \le c_0, 
$
where $c_0$ is some positive constant depending on $m$.
So,
$$
\int_{{\mathbb C^n} \times T_C} \frac {{\vert \Phi(z, \xi)\vert}^2 e^{-\varphi (z, \xi)}} {(1 + {\Vert (z, \xi)\Vert}^2)^{3n}} \  d\lambda_{2n} (z, \xi) \leq B e^{c_0}\int_{T_C} {\vert F(z)\vert}^2 e^{-\widetilde{\varphi}(z)} \ d \lambda_n (z).
$$
The right side of this inequality denote by $B_F$. 
For brevity put
$$
h_m(z, \xi) = 2(b(Im \xi)  + \omega_m (\Vert z \Vert) + m \ln (1 + \frac{1}{d(\xi)})), \ (z, \xi) \in {\mathbb C^n} \times T_C.
$$
Let us now obtain uniform estimates on $\Phi$ using plurisubharmonicity of  
$\vert \Phi(z, \xi)\vert$ in ${\mathbb C^n} \times T_C$. 
Let $(z, \xi) \in {\mathbb C^n \times T_C}$ and $R=\min (1, \frac {d(\xi)}{4})$.   
Note that if $(t, u) \in {\mathbb C^n \times T_C}$ belongs to the ball $B_R((z, \xi))$ then taking into account that:

1). for each $\varepsilon > 0$ there exists a constant $A_{\varepsilon}$ (by lemma 3) such that 
$$
\vert b(Im u) - b(Im \xi)\vert \le \varepsilon \Vert Im u \Vert + \varepsilon \Vert Im \xi \Vert + A_{\varepsilon}
$$
and, consequently, 
$$
\vert b(Im u) - b(Im \xi)\vert \le 2 \varepsilon \Vert Im \xi \Vert + A_{\varepsilon} + \varepsilon ;
$$

2). $$
\vert \ln(1 + {\Vert (t, u)\Vert}^2) - \ln(1 + {\Vert (z, \xi)\Vert}^2)\vert \le 1 ;
$$

3). 
$$
\left \vert \ln \left(1+\frac{1}{d(u)}\right) -  \ln \left (1+\frac{1}{d(\xi)}\right) \right \vert \le \frac 5 4.
$$

4). by inequality (3)
$$
\vert \omega_m (\Vert t \Vert) - \omega_m (\Vert z \Vert)\vert \le K_m,
$$
and putting $B_{\varepsilon}=2A_{\varepsilon} + 2\varepsilon +   \frac {5m} {2} + 2K_m + 5n + 2$,
we obtain for each $\varepsilon > 0$
$$
\vert 
h_m(t, u) + (5n + 2) \ln (1 + {\Vert (t, u)\Vert}^2) - h_m (z, \xi) - (5n + 2) \ln (1 + \Vert (z, \xi)\Vert^2)\vert \le 
$$
$$
\le
2 \varepsilon \Vert Im \xi \Vert + B_{\varepsilon}.
$$

Furthermore, for plurisubharmonic in $\mathbb C^n \times T_C$ function $ \vert \Phi(t, u) \vert^2 $ we have
$$
\vert \Phi(z, \xi)\vert^2 \leq 
\frac 
{1}{\nu_{2n}(R)}
\displaystyle \int \limits_{B_R((z, \zeta))}
\vert \Phi(t, u)\vert^2 \ d\lambda_{2n}(t, u) \leq  
$$
$$
\leq 
\frac 
{1}{\nu_{2n}(R)}
\int \limits_{\mathbb C^n\times T_C} \vert\Phi(t, u)\vert^2 e^{-(h_m(t, u) + (5n + 2) \ln(1+\Vert(t, u)\Vert^2))}  
 \ d\lambda_{2n}(t, u) \times
$$
$$
\times e^{\sup \limits_{(t, u) \in B_R((z,\zeta))} (h_m(t, u) + (5n + 2) \ln(1+\Vert(t, u)\Vert^2))}  \le 
$$
$$
\le 
\frac 
{B_F}{\nu_{2n}(R)}  e^{h_m (z, \xi) + (5n + 2) \ln (1 + \Vert (z, \xi)\Vert^2) + 2 \varepsilon \Vert Im \xi \Vert + B_{\varepsilon}}
$$
From here it follows that for each $\varepsilon > 0$ there is a constant $c_1 =c_1(\varepsilon) > 0$ such that for $(z, \xi) \in {\mathbb C^n} \times T_C$
$$
\vert \Phi(z, \xi)\vert \leq c_1
e^{b(Im \xi)  + \omega_m (\Vert z \Vert) + (m + 2n)\ln (1 + \frac{1}{d(\xi)}) + (5n + 2) \ln (1 + \Vert (z, \xi)\Vert) + 
\varepsilon \Vert Im \xi \Vert}. \eqno (9)
$$
Since $\Phi(z, \xi)$ is an entire function of $z$, then expanding $\Phi(z, \xi)$ in powers of $z$ we have
$$
\Phi(z, \xi) = \displaystyle \sum_{\vert\alpha\vert\geq 0} C_\alpha (\xi) z^\alpha,  \ \xi\in T_C,  \ z\in {\mathbb C^n}.
$$
By the Cauchy formula for arbitrary $\alpha\in {\mathbb Z_+^n, \ \ R>0}$
$$
C_\alpha (\xi) = 
\frac{1}{(2\pi i)^n} \int_{\vert z_1\vert =R}\ldots\int_{\vert z_n\vert =R} \frac {\Phi (z, \xi)} {z_1^{\alpha_1 +1} \ldots z_n^{\alpha_n +1}} \ dz_1\ldots dz_n .
$$
From here it follows that $C_\alpha \in H(T_C)$. 
Using (9) we have for $\xi \in T_C$
$$
\vert C_\alpha(\xi) \vert \leq 
\frac 
{c_1 ((1 + \sqrt{n} R)(1 + \Vert \xi \Vert))^{5n+2} e^{b(Im \xi) + \varepsilon \Vert Im \xi \Vert + \omega_m(\sqrt{n} R)}
(1 + \frac{1}{d(\xi)})^{m+2n}} {R^{\vert \alpha\vert}} \ .
$$
Using lemma 2 we find a constant $c_2 = c_2(\varepsilon)>0$ such that for each
$ R>0$ and $\xi \in T_C$
$$
\vert C_\alpha (\xi)\vert \leq 
c_2 \frac 
{e^{\omega_{m + 5n +2}(\sqrt{n} R)}} {R^{\vert\alpha\vert}} 
e^{b(Im \xi)+ \varepsilon \Vert Im \xi \Vert} (1 + \Vert\xi\Vert)^{5n+2} \left(1 + \frac{1}{d(\xi)}\right)^{m+2n}.
$$
Since 
$$
\displaystyle \inf_{R>0} \frac {e^{\omega(r)}} {r^k} = \frac{1} {M_k} \ ,  \ k = 0,1, \ldots,
$$
then for each $ \alpha\in {\mathbb Z_+^n}, \ 
\xi\in T_C$ and for each $\varepsilon > 0$ 
$$
\vert C_\alpha (\xi)\vert \leq 
c_2(\varepsilon) \left(\frac{\sqrt{n}} {\varepsilon_{m + 5n + 2}}\right)^{\vert\alpha\vert} 
\frac{e^{b(Im \xi)+ \varepsilon \Vert Im \xi \Vert}} {M_{\vert\alpha\vert}} 
(1 + \Vert\xi\Vert)^{5n+2} \left(1 + \frac{1}{d(\xi)}\right)^{m+2n}. 
$$
So, for each $\alpha\in {\mathbb Z_+^n}$ and for each $\varepsilon > 0$ \ $C_\alpha \in H_{b_{\varepsilon}}(T_C)$ and, consequently, $C_\alpha \in H_b(T_C)$.  According to the Corollary 3 for each $\alpha \in {\mathbb Z_+^n}$ \ $C_\alpha \in V_b(T_C)$.
By theorem 1 there exist functionals $S_\alpha \in S^*(U)$ such that $\hat {S_\alpha} = C_\alpha$. 

From the previous inequality also it follows that for each $\varepsilon > 0$ the set 
$\{M_{\vert\alpha\vert} \left(\frac {\varepsilon_{m + 5n + 2}} {\sqrt{n}} \right)^{\vert\alpha\vert} C_\alpha \}_{\alpha \in \mathbb Z_+^n} $ is bounded in $V_{b_{\varepsilon}}(T_C)$. Hence, it is bounded in $H_b(T_C)$. Since by Corollary 3 $H_b(T_C) = V_b(T_C)$, then this set is  bounded in $V_b(T_C)$.
In view of topological isomorphism between the spaces $S^*(U)$ and $V_b(T_C)$ the set 
${\cal A}= \{M_{\vert\alpha\vert} \left(\frac {\varepsilon_{m + 5n + 2}} 
{\sqrt{n}} \right)^{\vert\alpha\vert} S_\alpha \}_{\alpha \in \mathbb Z_+^n}$
is bounded in $S^* (U)$. Hence, it is weakly bounded.
By Schwartz theorem [1] there exist numbers $c_3 >0$ and $p\in\mathbb N$ such that
$$
\vert (F, f) \vert \leq c_3 {\Vert \varphi \Vert}_{p, U}, \ F \in {\cal A}, \ f \in S(U).  
$$ 
So, for all
$\alpha \in {\mathbb Z_+^n}, f \in S(U)$ 
$$
\vert (S_{\alpha}, f) \vert \leq c_3 \left(\frac{\sqrt{n}}{\varepsilon_{m + 5n + 2}}\right)^{\vert\alpha\vert} \frac{{\Vert  f \Vert}_{p, U}}{M_{\vert\alpha\vert}} \ . \eqno (10)
$$
Define the functional $T$ on $G_M(U)$ by the rule:
$$
(T, f) = \sum_{\vert\alpha\vert \geq  0} (S_\alpha, (-i)^{\vert\alpha\vert} D^\alpha f), \ f \in G_M(U).  \eqno (11)
$$
It is correctly defined. 
Using (10) we have for all $ f \in G_M(U), \alpha \in {\mathbb Z_+^n}, s \in {\mathbb N}$
$$
\vert (S_\alpha, D^\alpha f) \vert \leq 
c_3 \left(\frac {\sqrt{n}} {\varepsilon_{m + 5n + 2}}\right)^{\vert\alpha\vert} \frac{1} {M_{\vert\alpha\vert}}\sup_{x \in U, \vert\beta\vert\leq p }\vert (D^{\alpha+\beta} f)(x)\vert (1 + \Vert x\Vert)^p\leq 
$$
$$
\leq c_3 \left(\frac {\sqrt{n}} {\varepsilon_{m + 5n + 2}}\right)^{\vert\alpha\vert} \frac{1} {M_{\vert\alpha\vert}}\sup_{x\in U, \vert\beta\vert\leq p } \frac {p_s (f) \varepsilon_s^{\vert\alpha\vert +\vert\beta\vert} M_{\vert\alpha\vert +\vert\beta\vert}(1 + \Vert x\Vert)^p}  {(1 + \Vert x\Vert)^s} \ .
$$
It is easy to check that
$$
M_{k+s}\leq H_1^s H_2^{(k+s)s} M_k, \ k,s \in {\mathbb Z_+}.
$$
So for $s\geq p$ 
$$
\vert (S_\alpha, D^\alpha f) \vert \leq c_3 
\left(\frac {\sqrt{n}}{\varepsilon_{m + 5n + 2}}\right)^{\vert\alpha\vert} \frac{p_s (f)}{M_{\vert\alpha\vert}} 
\sup_{x\in U,\vert\beta\vert\leq p} \frac {\varepsilon_s^{\vert\alpha\vert + \vert\beta\vert} 
H_1^ {\vert\beta\vert} H_2^{(\vert\alpha\vert + \vert\beta\vert)\vert\beta\vert}M_{\vert\alpha\vert} } {(1 + \Vert x\Vert)^{s-p}} \le
$$
$$
\leq c_3 \left(\frac {\sqrt{n}}{\varepsilon_{m + 5n + 2}}\right)^{\vert\alpha\vert} H_1^p p_s (f) \varepsilon_s^{\vert\alpha\vert}h_2^{p^2}H_2^{p\vert\alpha\vert} 
= c_3 H_1^p H_2^{p^2} \left(\frac {\sqrt{n}\varepsilon_s H_2^p }  {\varepsilon_{m + 5n + 2}}\right)^{\vert\alpha\vert} p_s (f). 
$$
Choose $s$ so that
$\frac {\sqrt{n}\varepsilon_s H_2^p }  {\varepsilon_{m + 5n + 2}} = q_s <1$. 
Then for each $f \in G_M(U), \alpha \in {\mathbb Z_+^n}$ \ 
$$
\vert (S_\alpha, D^\alpha f) \vert \leq c_4 q_s^{\vert\alpha\vert} p_s (f), 
$$
where $c_4= c_3 H_1^p H_2^{p^2}$.
From here it follows that the series in (11) converges and 
$$
\vert (T, f) \vert\leq \frac {c_4 }{(1 - q)^n} p_s (f),\ \ f \in G_M(U). 
$$
Hence, the linear functional $T$ is correctly defined and continuous. 
Besides that $\hat T  = F.$ Indeed, 
$$
\hat  T(z)= \sum_{\vert\alpha\vert\leq 0} (S_{\alpha}, (-i)^{\vert\alpha\vert} D^{\alpha} (e^{i<\xi, z>})) 
=
\sum_{\vert\alpha\vert\leq 0}(S_{\alpha}, (-i)^{\vert\alpha\vert}(iz)^{\alpha}(e^{i<\xi, z>}))
$$
$$
=
\sum_{\vert\alpha\vert\leq 0}z^{\alpha}(S_{\alpha},(e^{i<\xi, z>}))=
\sum_{\vert\alpha\vert\leq 0} C_{\alpha}(z) z^{\alpha} = \Phi(z, z) = F(z).
$$
Thus, $L$ is surjective.

Now let us show that $L$ is injective. Let for $T \in G_M'(U)$ \ $\hat{T} \equiv~0$. We will show that $T$ is a zero functional. 
There exist numbers $m \in \mathbb N$ and $c>0$ such that
$$
\vert (T, f) \vert \leq c p_m (f), \ f \in G_M(U).  
$$
By lemma 4 there exist functionals $T_{\alpha}\in C_m'(U)$ ($\alpha\in\mathbb Z_+^n$) such that
$$
(T, f) = \displaystyle \sum_{\alpha \in \mathbb Z_+^n} (T_{\alpha}, D^{\alpha} f), \ f \in G_M(U),
$$ 
and 
$$
\vert (T_{\alpha}, g) \vert \leq \frac {c}  {\varepsilon_m^{\vert\alpha\vert} M_{\vert\alpha\vert}} \widetilde{p}_m (g), \ 
g \in C_m(U). \eqno (12)
$$
Hence, 
$$
\hat T(z) = \sum_{\alpha \in \mathbb Z_+^n} (T_{\alpha}, (iz)^{\alpha} e^{i<\xi, z>}) = 
\sum_{\alpha \in \mathbb Z_+^n} i^{\vert\alpha\vert} (T_{\alpha}, e^{i<\xi, z>}) z^{\alpha}, \ z \in T_C.
$$
Let $V_{\alpha}(z) = i^{\vert\alpha\vert} (T_{\alpha}, e^{i<\xi, z>}) $. 
Obviously, $V_{\alpha}\in H(T_C)$.
Using (12) and lemma~1, we obtain 
$$
\vert V_{\alpha} (z) \vert \leq \frac {{\mu}_m} {\varepsilon_m^{\vert\alpha\vert} M_{\vert\alpha\vert}}(1 + \Vert z \Vert)^{2m}
\left(1 + \frac {1}{\Delta_C(y)}\right)^{3m} e^{b(y)}, \eqno (13)
$$
where  ${\mu}_m >0$ is some constant independent of $z = x + iy  \in T_C$. 
Let
$$
S (u, z) = \displaystyle \sum_{\vert\alpha\vert\geq 0} V_{\alpha} (z) u^{\alpha}, \ z \in T_C, u \in \mathbb C^n.  
$$
Using (13) we obtain
$$
\vert S(u, z)\vert  
\leq   \sum_{\vert \alpha\vert \geq 0} \frac {{\mu}_m \Vert u\Vert^{\vert\alpha\vert}}  {\varepsilon_m^{\vert\alpha\vert} M_{\vert\alpha\vert}} (1 + \Vert z \Vert)^{2m} \left(1 + \frac {1}{\Delta_C(y)}\right)^{3m} e^{b(y)}=  
$$
$$
= {\mu}_m \left(1 + \frac {1}{\Delta_C(y)}\right)^{3m} e^{b(y)}
\sum_{\vert \alpha\vert \geq 0}\frac {\Vert u\Vert^{\vert\alpha\vert}}  {\varepsilon_{m+1}^{\vert\alpha\vert} M_{\vert\alpha\vert}} \left(\frac {\varepsilon_{m+1}}  {\varepsilon_m}\right)^{\vert\alpha\vert} \leq
$$
$$
\leq {\mu}_m e^{b(y)}\left(1 + \frac {1}{\Delta_C(y)}\right)^{3m}
\sup_{\alpha\in\mathbb Z_+^n} \frac {\Vert u\Vert^{\vert\alpha\vert}}  {\varepsilon_{m+1}^{\vert\alpha\vert} M_{\vert\alpha\vert}} \sum_{\vert\alpha\vert\geq 0} \left(\frac {\varepsilon_{m+1}}  {\varepsilon_m}\right)^{\vert\alpha\vert} =  
$$
$$
\leq {\mu}_m e^{b(y)}\left(1 + \frac {1}{\Delta_C(y)}\right)^{3m}
e^{\omega_{m+1}(\Vert u\Vert)} \left(\frac {\varepsilon_m}  {\varepsilon_m - \varepsilon_{m+1}}\right)^n \ .
$$
Note that
$$
S(z, z) = \sum_{\vert\alpha\vert\geq 0} V_{\alpha}(z) z^{\alpha} = 0, \ \forall z\in T_C.
$$ 
Then by theorem 2 there exist functions 
$S_1,\ldots,S_n \in H(\mathbb C^n \times T_C)$  and a number $p \in {\mathbb N}$ such that
$$
S(z, \zeta) = \sum_{j=1}^n S_j(z, \zeta) (z_j-\zeta_j), \ z \in \mathbb C^n, \zeta \in T_C,
$$ 
and for each $\varepsilon > 0$ there exists a number $a_{\varepsilon} > 0 $ such that for $j=1, \ldots , n$ and
for $(z, \zeta) \in {\mathbb C^n \times T_C}$
$$
\vert S_j(z, \zeta) \vert \leq a_{\varepsilon} \exp ( \omega_p(\Vert z\Vert) + b(Im \zeta)+ \varepsilon \Vert Im \zeta \Vert  + p \ln (1 + \frac{1} {d(\zeta)}) + p\ln (1 + \Vert\zeta\Vert)). 
$$
Further, we expand  $S_j$ in powers of $z$:
$$
S_j(z, \zeta) = \sum_{\vert\alpha\vert\geq 0} S_{j, \alpha} (\zeta) z^{\alpha},\ \ z\in\mathbb C^n, \ \zeta\in T_C.  
$$
From the last inequality we have 
$$
\vert S_{j, \alpha}(\zeta) \vert 
\leq \inf_{R>0}\frac {\displaystyle\max_{\vert z_1\vert=R,\ldots,\vert z_n\vert=R} \vert S_j(z, \zeta)\vert}  {R^{\vert\alpha\vert}}\leq 
$$
$$
\leq \inf_{R>0} \frac {a_{\varepsilon}\exp(\omega_p(\sqrt{n}R)) e^{(b(Im \zeta) + \varepsilon \Vert Im \zeta \Vert + p\ln(1+\frac{1}{d(\zeta)}) + p\ln (1+\Vert\zeta\Vert))}}  {R^{\vert\alpha\vert}} = 
$$
$$
= a_{\varepsilon}\exp(b(Im\zeta) + \varepsilon \Vert Im \zeta \Vert + p\ln \left(1+\frac{1}{d(\zeta)}\right) + p\ln (1+\Vert\zeta\Vert)) \inf_{R>0}\frac {e^{\omega (R)} \sqrt{n}^{\vert\alpha\vert}}  {(R\varepsilon_p)^{\vert\alpha\vert}} = 
$$
$$
= a_{\varepsilon} \exp((b(Im\zeta) +\varepsilon \Vert Im \zeta \Vert  + p\ln \left(1+\frac{1}{d(\zeta)}\right) + p\ln (1+\Vert\zeta\Vert)) \frac {1} {M_{\vert\alpha\vert}} 
\left(\frac {\sqrt{n}} {\varepsilon_p}\right)^{\vert\alpha\vert} \ .   
$$
Choose $k\in\mathbb N$ so that $\varepsilon_k\sqrt{n} < \varepsilon_p$. 
Then for $\zeta\in T_C$
$$
\vert S_{j,\alpha}(\zeta)\vert \leq \frac {a_{\varepsilon} \exp(b(Im \zeta) + \varepsilon \Vert Im \zeta \Vert  + p\ln(1+\frac{1}{d(\zeta)}) + p\ln (1+\Vert\zeta\Vert))}  {\varepsilon_k^{\vert\alpha\vert} M_{\vert\alpha\vert}}.
$$
Thus, for all $\alpha\in {\mathbb Z_+^n}, j = 1, \ldots , n$ and for each $\varepsilon > 0$ \ $S_{j,\alpha} \in V_{b_{\varepsilon}}(T_{C})$. According to the corollary 3 of the theorems 1 $S_{j,\alpha} \in V_b(T_C)$ for each $\alpha\in {\mathbb Z_+^n}, j = 1, \ldots , n$.
Since the Fourier-Laplace transform establishes a topological isomorphism between the spaces $S^*(U)$ and $V_b(T_C)$, 
then there exist functionals $\psi_{j, \alpha} \in S^*(U)$ such that $\hat \psi_{j, \alpha} = S_{j, \alpha}$. 
From the last inequality it follows that the set 
$\{S_{j, \alpha}\varepsilon_k^{\vert\alpha\vert} M_{\vert\alpha\vert} \}_{\alpha \in \mathbb Z_+^n} $ is bounded in each space 
$V_{b_{\varepsilon}}(T_C)$. But then it is bounded in $H_b(T_C)$. Now by  Corollary 3 of the theorems 1 this set is bounded in $V_b(T_C)$.
Then the set 
$
\Psi = \{\varepsilon_k^{\vert\alpha\vert} M_{\vert\alpha\vert}\psi_{j, \alpha} \}_{\alpha \in \mathbb Z_+^n, j=1,\ldots,n}
$ 
is bounded in $S^*(U)$. And then it is weakly bounded in $S^*(U)$.
By Schwartz theorem [1]  there exist numbers $c_4>0$ and $p\in\mathbb N$ such that
$$
\vert (F, \varphi) \vert \leq c_4 \vert \Vert \varphi \Vert\vert_{p, U}, \ F \in \Psi, \ \varphi\in S(U).  
$$
Thus, $\forall j=1,\ldots,n,\ \ \forall \alpha\in\mathbb Z_+^n$ 
$$
\vert (\Psi_{j, \alpha}, f) \vert\leq \frac{c_4} {\varepsilon_k^{\vert\alpha\vert} M_{\vert\alpha\vert}} {\Vert f\Vert}_{p, U}, \ 
f \in S(U). \eqno (14)
$$ 
For $j=1,\ldots,n$ and $\alpha \in \mathbb Z^n$ with at least one negative component 
let $\Psi_{j, \alpha}$ be a zero functional from $S^*(U)$ and $S_{j, \alpha} (z)=0,\ \ \forall z \in \mathbb C^n$. 
Then
$$
S(z, \zeta) = \sum_{j=1}^n S_j(z, \zeta) (z_j - \zeta_j) = 
\sum_{j=1}^n \sum_{\vert\alpha\vert\geq 0} S_{j, \alpha}(\zeta) z^{\alpha} (z_j - \zeta_j)=
$$ 
$$
=\sum_{j=1}^n \sum_{\vert\alpha\vert\geq 0} S_{j,\alpha}(\zeta)z_1^{\alpha_1}\ldots z_j^{\alpha_j +1} \ldots z^{\alpha_n}-\sum_{j=1}^n \sum_{\vert\alpha\vert\geq 0}S_{j,\alpha}(\zeta)z^{\alpha}\zeta_j=  $$ 
$$
=\sum_{j=1}^n \sum_{\vert\alpha\vert\geq 0} (S_{j,(\alpha_1,\ldots,\alpha_{j}-1,\ldots,\alpha_n)}(\zeta) - S_{j, \alpha}(\zeta) \zeta_j)z^{\alpha},\ \ z\in\mathbb C^n, \ \zeta\in T_C. 
$$ 
Since 
$$
S(z, \zeta) = \sum_{\vert\alpha\vert\geq 0} V_{\alpha}(\zeta)z^{\alpha},\ \ \forall \zeta\in T_C,  \ z\in\mathbb C^n,
$$ we have $\forall\alpha\in\mathbb Z_+^n$
$$V_{\alpha}(\zeta) = \sum_{j=1}^n (S_{j,(\alpha_1,\ldots,\alpha_{j-1},\ldots,\alpha_n)} (\zeta) - S_{j, \alpha}(\zeta)\zeta_j). 
$$ 
The expression on the right can be represented as
$$
\sum_{j=1}^n (\hat \Psi _{j,(\alpha_1,\ldots,\alpha_{j-1},\ldots,\alpha_n)}(\zeta) + 
i (\Psi_{j, \alpha}, \frac {\partial} {\partial \xi_j}(e^{i<\xi, \zeta>})) ).
$$  
That is, the right side is the Fourier-Laplace of a functional acting by the rule
$$
f \in S(U) \rightarrow \sum_{j=1}^n(\Psi_{j, (\alpha_1,\ldots, \alpha_{j-1}, \ldots, \alpha_n)}, f) + 
i(\Psi_{j, \alpha}, (\frac{\partial} {\partial\xi_j}f)).  
$$
This means that 
$$
(T_{\alpha}, f) = (-i)^{\vert\alpha\vert} \sum_{j=1}^n (i(\Psi_{j, \alpha}, (\frac {\partial} {\partial\xi_j} f)) + (\Psi_{j, (\alpha_1, \ldots, \alpha_{j-1},\ldots, \alpha_n)}, f )).
$$
Thus, for $f \in G_M(U)$
$$
(T, f) = \sum_{\vert\alpha\vert\geq 0} (T_{\alpha}, D^{\alpha} f) = 
$$
$$
=\sum_{\vert\alpha\vert\geq 0}(-i)^{\vert\alpha\vert} \sum_{j=1}^n (i(\Psi_{j, \alpha}, (\frac {\partial} {\partial\xi_j} D^{\alpha}f)) + (\Psi_{j,(\alpha_1,\ldots,\alpha_{j-1},\ldots,\alpha_n)},D^{\alpha} f)), 
$$
For arbitrary $N \in \mathbb N$ and $j =1, \ldots , n)$
let us define the sets 
$$
B_N = \{\alpha=(\alpha_1,\ldots,\alpha_n)\in \mathbb Z^n: \alpha_1\leq N,\ldots,\alpha_n\leq N \},
$$
$$R_{N, j} = \{\alpha_1\leq N,\ldots,\alpha_j = N,\ldots,\alpha_n\leq N, \alpha\in\mathbb Z_+^n\}$$ and a functional $T_N$ on $G_M(U)$ by the rule
$$
(T_N, f) = \sum_{\alpha \in B_N} (-i)^{\vert\alpha\vert} 
\sum_{j=1}^n (i(\Psi_{j, \alpha}, (\frac {\partial} {\partial\xi_j} D^{\alpha}f)(\xi)) + (\Psi_{j,(\alpha_1,\ldots,\alpha_{j-1},\ldots,\alpha_n)},D^{\alpha} f)).
$$
Then $(T, f) = \displaystyle\lim_{N\rightarrow\infty} (T_N, f),\ \ f\in G_M(U)$. 

From the representation
$$
(T_N, f) = 
\sum_{j=1}^n (\sum_{\alpha\in B_N} ((-i)^{\vert\alpha\vert} i (\Psi_{j, \alpha},\frac{\partial} {\partial\xi_j} D^{\alpha} f) + 
$$
$$
+
\sum_{\beta\in B_N}(-i)^{\vert\beta\vert} (\Psi_{j,(\beta_1,\ldots,\beta_{j-1},\ldots,\beta_n)}, D^{\beta} f)). 
$$
we obtain that for fixed $j \in \{1, \ldots, n \}$ terms corresponding to the multi-index $\alpha$ with
$\alpha_1 \leq N, \ldots,\alpha_j \leq N-1, \ldots,\alpha_n \leq N$ and the terms corresponding to he multi-index 
$\beta = (\beta_1, \ldots, \ldots, \beta_n)= \alpha_n$ with $\beta_1 = \alpha_1, \ldots, \beta_j = \alpha_j+1, \ldots, \beta_n= \alpha_n$ mutually vanishes each other. 
From this we have
$$
(T_N, f) = \sum_{j=1}^n \sum_{\alpha\in R_{N, j}}(-i)^{\vert\alpha\vert} 
i (\Psi_{j, \alpha},\frac {\partial} {\partial \xi_j} D^{\alpha}f), \ f \in G_M(U).
$$
Further, taking into account (14), $\forall f\in G_M(U)$ 
$$
\vert (T_N, f) \vert \leq\sum_{j=1}^n \sum_{\alpha \in R_{N, j}} \vert (\Psi_{j, \alpha},\frac {\partial} {\partial \xi_j} D^{\alpha}f)  \vert 
\leq   
$$
$$
\leq \sum_{j=1}^n \sum_{\alpha \in R_{N, j}} \frac {c_4} {\varepsilon_k^{\vert\alpha\vert} M_{\vert\alpha\vert}} 
\sup_{\xi\in U,\vert\gamma\vert\leq p} (\vert D^{\gamma}(\frac {\partial} {\partial \xi_j}D^{\alpha} f)(\xi) \vert (1 + \Vert \xi \Vert)^p)= 
$$ 
$$
= \sum_{j=1}^n \sum_{\alpha \in R_{N, j}} 
\frac {c_4} {\varepsilon_k^{\vert\alpha\vert} M_{\vert\alpha\vert}} \sup_{\xi\in U, \vert \gamma \vert \leq p}(\vert(D^{(\alpha_1+\gamma_1, \ldots\alpha_j+\gamma_j+1,\ldots,\alpha_n + \gamma_n)} f)(\xi)\vert (1 + \Vert \xi \Vert)^p)).
$$ 
Choose natural $s >k$ so that $q=\frac{\varepsilon_sH_2^{p+1}} {\varepsilon_k}  < 1 $. 
Then $\forall f \in G_M(U)$ 
$$
\vert (T_N, f) \vert \leq \sum_{j=1}^n \sum_{\alpha \in R_{N, j}} \frac {c_4} {\varepsilon_k^{\vert\alpha\vert} M_{\vert\alpha\vert}} p_s (f)
\sup_{\xi \in U, \vert\gamma\vert\leq p} \frac 
{\varepsilon_s^{\vert\alpha\vert+\vert\gamma\vert+1} M_{\vert\alpha\vert+\vert\gamma\vert+1}}  {(1 + \Vert\xi\Vert)^{s-p}} \leq 
$$ 
$$
\leq \sum_{j=1}^n \sum_{\alpha \in R_{N, j}} \frac {c_4} {\varepsilon_k^{\vert\alpha\vert} M_{\vert\alpha\vert}} p_s (f) \varepsilon_s^{\vert\alpha\vert}\sup_{\vert\gamma\vert\leq p}H_1^{\vert\gamma\vert +1} H_2^{(\vert\alpha\vert+\vert\gamma\vert +1)(\vert\gamma\vert + 1)}M_{\vert\alpha\vert} \leq
$$ 
$$
\leq \sum_{j=1}^n \sum_{\alpha \in R_{N, j}} c_4 p_s (f) \left(\frac{\varepsilon_s} {\varepsilon_k} H_2^{p+1}\right)^{\vert\alpha\vert} H_1^{p+1} H_2^{(p+1)^2}<   
$$ 
$$
< c_4 H_1^{p+1} H_2^{(p+1)^2} p_s (f) n q^N (N+1)^{n-1}.
$$
From here it follows that $(T_N, f) \rightarrow 0$ at $N \rightarrow \infty, \ \forall f\in G_M(U)$. Hence $(T, f) =0, \ \forall f\in G_M(U)$. 
So, $T$ is a zero functional. 

By the open mapping  [8], [9] $L^{-1}$ is continuous. 
Thus,  $L$ is a topological isomorphism. 

{\bf Remark 2}. In the situation considered by Roever from his results [2, theorems 2.21.ii, 2.24.ii] it follows that $G_M^*(U)$ is topologically isomorphic to the projective limit of spaces $R_{C_{\varepsilon},
\varepsilon}$ ($\{C_{\varepsilon}\}$ is a family of compact subcones in cone $C$, $\varepsilon > 0)$, where $R_{C_{\varepsilon}, \varepsilon}$ is an inductive limit of spaces
$$
R_{C_{\varepsilon}, \varepsilon}^{(m)} = \{F \in H(T_C): {\Vert F \Vert}_{C_{\varepsilon},
\varepsilon}^{(m)} = \sup_{z \in T_{C_1}, \Vert y \Vert \ge
\varepsilon } \frac {\vert F(z) \vert} {e^{b(y) + \omega_m (\Vert
z \Vert)}} < \infty \} \ , \ m \in {\mathbb N}.
$$

\pagebreak

\vspace{0,5cm} 
{\sc Institute of mathematics with computer center of Ufa scientific center of Russian Academy of sciences, Ufa, RUSSIA}

{\it E-mail address}: musin@matem.anrb.ru

{\it E-mail address}: polina81@rambler.ru

\end{document}